# A Database of 2,500 Quasicrystal Cells

Authors:

Tony Robbin: wrote a program in the late 1980s implementing the deBruijn Dual Method for generating quasicrystals and used it to make a large sculpture in Denmark. (Photo 1, reference list)

George Frances: rewrote Robbin's program with improvements in the early 2000s and used it to generate an interactive, immersive display in a *Cave*. (Photos 2,3)

Kurt Bauman: this year got Robbin's 1980s computer to work and used it to generate a comma-delineated database of 2,500 quasicrystal cells that is presented here.

In the database, the first three numbers are the x, y, z components of the first vertex of the first cell, and so on for all eight sets of coordinates of the cell. The last number is a 1 or a 2 which indicates whether the cell is a thin or fat cell respectively. The program that generates the database computes the volume of each cell by its determinant, yielding only one of two numbers (0.47023 or 0.76084), numbers that are in the Golden Ratio. (A property of one-dimensional quasicrystals is that the length of the two unit cells --- line segments—are in the Golden ratio; a property of two-dimensional quasicrystals is that the area of the two unit cells –rhombs—are in the Golden ratio; and a property of three-dimensional quasicrystals is that the volumes of the two unit cells, parallelepipeds, are also in the Golden ratio). It was reassuring that calculations of each volume are true to 5 decimal places (see Appendix 1), but the authors felt that an integer value for the indicator would be easier for machines to read.

| | | | | | | |
|---|---|---|---|---|---|---|
| 1 | 4.471809E-01 | 1.376398E+00 | -1.036656E+01 | -4.472463E-01 | 1.376398E+00 | -1.081378E+01 |
| | 1.707902E-01 | 5.257465E-01 | -1.081378E+01 | 1.170791E+00 | 8.506713E-01 | -1.081378E+01 |
| | -7.236370E-01 | 5.257465E-01 | -1.126099E+01 | 2.763636E-01 | 8.506713E-01 | -1.126099E+01 |
| | 8.944001E-01 | 1.966953E-05 | -1.126099E+01 | -2.709031E-05 | 1.966953E-05 | -1.170820E+01 |
| | 1 | | | | | |
| 2 | -1.447230E+00 | 2.503395E-06 | -7.577709E+00 | -2.341657E+00 | 2.503395E-06 | -8.024922E+00 |
| | -1.723621E+00 | -8.506491E-01 | -8.024922E+00 | -7.236202E-01 | -5.257244E-01 | -8.024922E+00 |
| | -2.618048E+00 | -8.506491E-01 | -8.472136E+00 | -1.618047E+00 | -5.257244E-01 | -8.472136E+00 |
| | -1.000011E+00 | -1.376376E+00 | -8.472136E+00 | -1.894438E+00 | -1.376376E+00 | -8.919350E+00 |
| | 1 | | | | | |

The program that computes the database allows for different initial conditions (minute differences in starting positions), which generates a unique quasicrystal. For this quasicrystal, the 6 generating axes were offset from the origin by 0.1, 0.02, 01, 0.03, 0.04, 0.05, respectively. (see Appendix 1) These offsets generated a database in which there are an equal number of thin and fat unit cells and no sub-assemblies, no clusters that form rhombic polyhedrals such as a rhombic dodecahedron, or a rhombic icosahedron, or a rhombic triacontahedron. Another set of offsets would have a different outcome. (see Photos 3 and 5)

In every cell, the connections of the vertices are as follows: (see Photo 4)

1-2
1-4
1-6
2-3
2-5
3-4
3-8
4-7
5-6
5-8
6-7
7-8

The full database is linked here:

[tonyrobbin.net/pdfs/QuasicrystalCells2500.txt](tonyrobbin.net/pdfs/QuasicrystalCells2500.txt)   (Text file)

[tonyrobbin.net/pdfs/QuasicrystalCells2500.csv](tonyrobbin.net/pdfs/QuasicrystalCells2500.csv)   (csv file)

The database may also be useful for mathematicians and physicists for the further study of quasicrystals. (see Appendix 2)

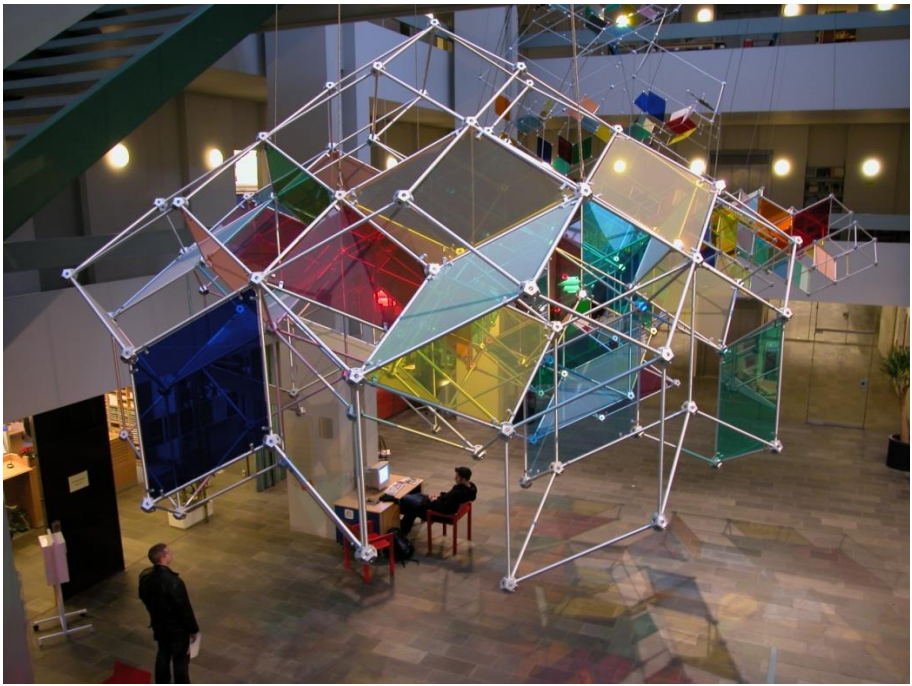

Photo 1

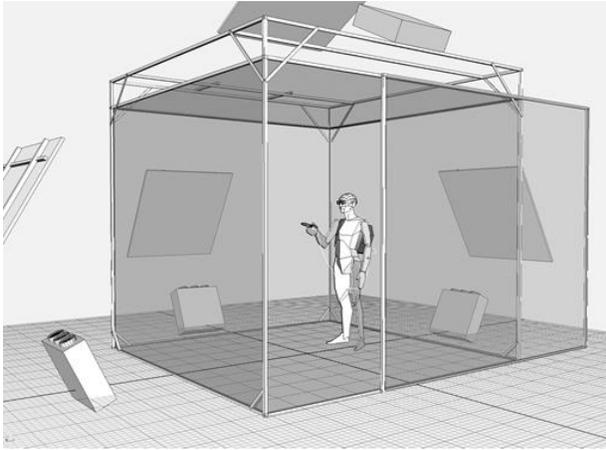
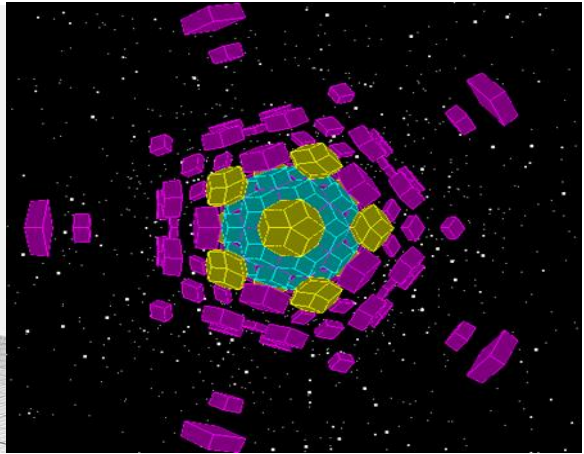

Photo 2                                                      Photo 3

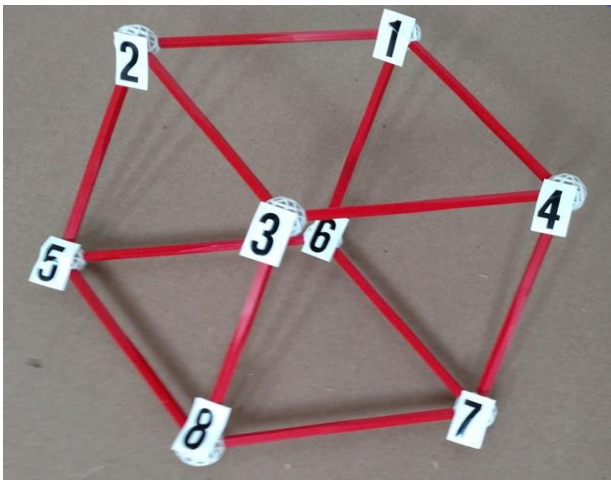
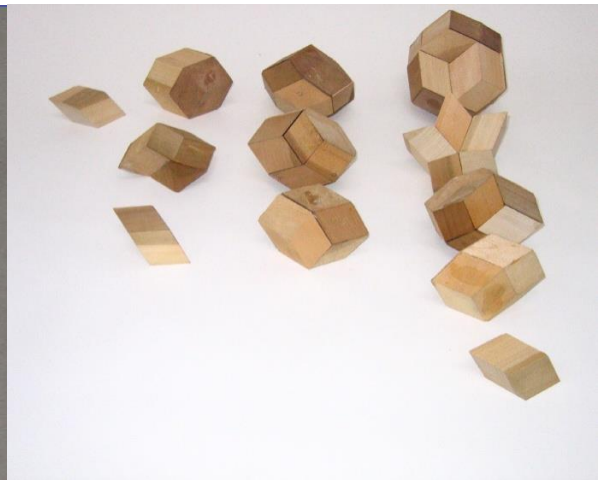

Photo 4 shows a cell with the vertices numbered.   Photo 5 shows the decomposition of rhombic polyhedra into constituent fat and thin cells.

The authors hope there will be a visual representation of the database in an open-access environment that would be sliceable (as flat slabs, domes, or arbitrary shapes) so that quasicrystal geometry could be used by sculptors and architects to make structures. Indeed, George Francis and his students present here a JavaScript version that will present a visual representation of the database on just about any computer, tablet, or smartphone. http://new.math.uiuc.edu/Tavares/htr4.html

Remaining projects are to convert this database to a format that could be directly imported into AutoCad or follow-on programs such as Rhino, Vectorworks, and the like. Additionally, the JavaScript program could be modified to recreate the functionality of the original Pascal to adjust the starting parameters, which will generate multiple distinct quasicrystal databases. The authors invite collaboration to realize these goals.

Appendix 1:

Sample raw data from the first 5 cells

1, 4.471809000E-01, 1.376398000E+00,-1.036656000E+01,-4.472463000E-01, 1.376398000E+00,-1.081378000E+01, 1.707902000E-01, 5.257465000E-01,-1.081378000E+01, 1.170791000E+00, 8.506713000E-01,-1.081378000E+01,-7.236370000E-01, 5.257465000E-01,-1.126099000E+01, 2.763636000E-01, 8.506713000E-01,-1.126099000E+01, 8.944001000E-01, 1.966953000E-05,-1.126099000E+01,-2.709031000E-05, 1.966953000E-05,-1.170820000E+01,-==4.70231076E-001==

2,-1.447230000E+00, 2.503395000E-06,-7.577709000E+00,-2.341657000E+00, 2.503395000E-06,-8.024922000E+00,-1.723621000E+00,-8.506491000E-01,-8.024922000E+00,-7.236202000E-01,-5.257244000E-01,-8.024922000E+00,-2.618048000E+00,-8.506491000E-01,-8.472136000E+00,-1.618047000E+00,-5.257244000E-01,-8.472136000E+00,-1.000011000E+00,-1.376376000E+00,-8.472136000E+00,-1.894438000E+00,-1.376376000E+00,-8.919350000E+00,-==4.70230577E-001==

3,-3.341641000E+00,-1.376393000E+00,-3.788854000E+00,-4.236068000E+00,-1.376393000E+00,-4.236068000E+00,-3.618032000E+00,-2.227045000E+00,-4.236068000E+00,-2.618031000E+00,-1.902120000E+00,-4.236068000E+00,-4.512459000E+00,-2.227045000E+00,-4.683281000E+00,-3.512458000E+00,-1.902120000E+00,-4.683281000E+00,-2.894422000E+00,-2.752772000E+00,-4.683281000E+00,-3.788849000E+00,-2.752772000E+00,-5.130495000E+00,-==4.70230391E-001==

4,-4.512450000E+00,-2.227051000E+00,-4.472136000E-01,-5.406877000E+00,-2.227051000E+00,-8.944272000E-01,-4.788841000E+00,-3.077703000E+00,-8.944272000E-01,-3.788840000E+00,-2.752778000E+00,-8.944272000E-01,-5.683268000E+00,-3.077703000E+00,-1.341641000E+00,-4.683267000E+00,-2.752778000E+00,-1.341641000E+00,-4.065231000E+00,-3.603430000E+00,-1.341641000E+00,-4.959658000E+00,-3.603430000E+00,-1.788854000E+00,-==4.70230668E-001==

5,-6.406861000E+00,-3.603447000E+00, 2.341641000E+00,-7.301288000E+00,-3.603447000E+00, 1.894427000E+00,-6.683252000E+00,-4.454099000E+00, 1.894427000E+00,-5.683251000E+00,-4.129174000E+00, 1.894427000E+00,-7.577679000E+00,-4.454099000E+00, 1.447214000E+00,-6.577678000E+00,-4.129174000E+00, 1.447214000E+00,-5.959642000E+00,-4.979825000E+00, 1.447214000E+00,-6.854069000E+00,-4.979825000E+00, 1.000000000E+00,-==4.70230515E-001==

And the last 5 cells

2495,-2.170802000E+00,-2.227047000E+00, 3.341641000E+00,-1.447200000E+00,-1.701310000E+00, 2.894427000E+00,-2.447205000E+00,-1.376400000E+00, 2.894427000E+00,-2.170802000E+00,-2.227047000E+00, 2.341641000E+00,-1.723603000E+00,-8.506622000E-01, 2.447214000E+00,-1.447200000E+00,-1.701310000E+00, 1.894427000E+00,-

2.447205000E+00,-1.376400000E+00, 1.894427000E+00,-1.723603000E+00,-8.506622000E-01, 1.447214000E+00,-7.60845739E-001

2496,-1.894376000E+00,-8.583209000E+00,-4.683281000E+00,-1.170774000E+00,-8.057471000E+00,-5.130495000E+00,-2.170779000E+00,-7.732561000E+00,-5.130495000E+00,-1.894376000E+00,-8.583209000E+00,-5.683281000E+00,-1.447177000E+00,-7.206824000E+00,-5.577709000E+00,-1.170774000E+00,-8.057471000E+00,-6.130495000E+00,-2.170779000E+00,-7.732561000E+00,-6.130495000E+00,-1.447177000E+00,-7.206824000E+00,-6.577709000E+00,-7.60846284E-001

2497,-2.341595000E+00,-7.206830000E+00,-2.788854000E+00,-1.617993000E+00,-6.681093000E+00,-3.236068000E+00,-2.617999000E+00,-6.356183000E+00,-3.236068000E+00,-2.341595000E+00,-7.206830000E+00,-3.788854000E+00,-1.894396000E+00,-5.830445000E+00,-3.683281000E+00,-1.617993000E+00,-6.681093000E+00,-4.236068000E+00,-2.617999000E+00,-6.356183000E+00,-4.236068000E+00,-1.894396000E+00,-5.830445000E+00,-4.683281000E+00,-7.60845807E-001

2498,-1.894387000E+00,-5.830452000E+00,-4.472136000E-01,-1.170785000E+00,-5.304714000E+00,-8.944272000E-01,-2.170791000E+00,-4.979804000E+00,-8.944272000E-01,-1.894387000E+00,-5.830452000E+00,-1.447214000E+00,-1.447188000E+00,-4.454067000E+00,-1.341641000E+00,-1.170785000E+00,-5.304714000E+00,-1.894427000E+00,-2.170791000E+00,-4.979804000E+00,-1.894427000E+00,-1.447188000E+00,-4.454067000E+00,-2.341641000E+00,-7.60846392E-001

2499,-7.235695000E-01,-4.979800000E+00, 1.447214000E+00, 3.266335000E-05,-4.454062000E+00, 1.000000000E+00,-9.999727000E-01,-4.129153000E+00, 1.000000000E+00,-7.235695000E-01,-4.979800000E+00, 4.472136000E-01,-2.763705000E-01,-3.603415000E+00, 5.527864000E-01, 3.266335000E-05,-4.454062000E+00, 0.000000000E+00,-9.999727000E-01,-4.129153000E+00, 0.000000000E+00,-2.763705000E-01,-3.603415000E+00,-4.472136000E-01,-7.60845985E-001

2500,-1.170789000E+00,-3.603422000E+00, 3.341641000E+00,-4.471865000E-01,-3.077684000E+00, 2.894427000E+00,-1.447192000E+00,-2.752774000E+00, 2.894427000E+00,-1.170789000E+00,-3.603422000E+00, 2.341641000E+00,-7.235897000E-01,-2.227037000E+00, 2.447214000E+00,-4.471865000E-01,-3.077684000E+00, 1.894427000E+00,-1.447192000E+00,-2.752774000E+00, 1.894427000E+00,-7.235897000E-01,-2.227037000E+00, 1.447214000E+00,-7.60845848E-001

Additional output from the program's counters for this set of offsets:

2500 number of cells

0 number of singularities (by which is meant locations where more than three planes intersect, that is locations where cells cluster to make rhombic polyhedrals.

1250 number of FAT cells

1250 number of THIN cells

Appendix 2: The construction of the database and what it means:

Imagine a 6-dimensional grid marked off from -25 to +25 on each of the six axes. Construct planes perpendicular to these axes at all the points marked off. Now imagine this structure projected to 3 dimensions such that instead of the axes being mutually perpendicular, the 90-degree angles are now (about) 63.44 degrees, meaning that the origin of the structure is at the center of a dodecahedron, the axes pass through the centers of the faces of the dodecahedron, and the planes are all parallel to the faces. Every two planes intersect in a line and a third plane intersects that line at a point, call it *intersect*. Use that point to construct a parallelepiped; say, use the point as the lower-left corner and grab the surrounding plane sections to make a skewed box (let's now call them cells). If you do this for all the combinations of axes and planes, the cells will be massively interlaced. Nicolaas deBruijn saw that it is possible to select only those cells that tesselate perfectly in three dimensions by referring back to the original pre-projected grid structure. Take a vector from the origin inside the dodecahedron to *intersect* and lay it along the other three axes that pass through the faces of the dodecahedron, the ones not used to calculate *intersect* in the first place. Take note of those positions on the axes. Now you have six integers referring to six planes of the original grid structure, making coordinates that identify a unique 3-dimensional cell in 6-dimensional space. If you project and plot this cell and only others constructed by this method, a perfect quasicrystal will be formed. (Details and complications abound as you program this.)

We usually think of projection as a one-way gate and that information is lost by the act of projection. The deBruijn method shows that this is not necessarily the case, and that information from the pre-projection state can pass through the gate, or one could say the algorithm passes back and forth between the two states, or if this were modeling a physical system, one could say that time flows both forward and backward in the algorithm.

Reconsider the 6-dimensional rectilinear grid first constructed. Visualize the 23rd box along the fifth dimension. You know its size, shape, orientation, and how it is connected to its neighbors. You know these far off details standing at the origin. These are exactly the details of far off (skewed) boxes that are so mysterious in quasicrystals. If nature uses something like the deBruijn algorithm to make experience, then quantum, non-local phenomena are not so mysterious. In other words, the method could be a useful model for such phenomena: it shows how there can be information at a distance where there is seemingly no continuity of space. (Tony Robbin)